\documentclass[11pt]{article}
\newcommand{\tabincell}[2]{\begin{tabular}{@{}#1@{}}#2\end{tabular}}

\date{}
\makeatletter
\@addtoreset{equation}{section}
\makeatother
\makeatother
\usepackage{amssymb,amsmath}
\usepackage{color}
\usepackage[colorlinks,linkcolor=blue,citecolor=blue]{hyperref}
\usepackage{cite}
\usepackage{graphicx}
\usepackage{multirow}
\usepackage{array}
\usepackage{booktabs} 
\usepackage{times}

 \numberwithin{equation}{section}
 \newtheorem{thm}{Theorem}[section]
\newtheorem{lem}[thm]{Lemma}
\newtheorem{rem}[thm]{Remark}

\allowdisplaybreaks

\begin{document}
\title{\LARGE  Ground states for the double weighted critical Kirchhoff equation on the unit ball in $\mathbb{R}^3$}
\author{Yao Du$^a$,\ \ \ Jiabao Su$^b$\footnote{Corresponding author. \hfill\break\indent \ \  E-mail addresses: 1052591976@qq.com(Y. Du),\ sujb@cnu.edu.cn(J. Su). }\\ {\small  $^a$School of Science, Xihua University} \\
{\small Chengdu 610039, People's Republic of China}\\
{\small  $^b$School of Mathematical Sciences, Capital Normal University}\\ {\small Beijing 100048, People's Republic of China}}
\maketitle

\begin{abstract}
 This paper deals with the existence of ground states for degenerative ($a=0$) and non-degenerative ($a>0$) double weighted critical Kirchhoff equation
\begin{eqnarray*}
\left\{
  \begin{array}{ll}
\displaystyle-\left(a+b\int_B |\nabla u|^2dx\right)\Delta u=|x|^{\alpha_1} |u|^{4+2\alpha_1}u+\mu|x|^{\alpha_2} |u|^{4+2\alpha_2}u+\lambda h(|x|) f(u)  &{\rm in}\ B,\\
u=0  &{\rm on}\ \partial B,
 \end{array}
\right.
\end{eqnarray*}
where $B$ is a unit open ball in $\mathbb{R}^3$ with center $0$, $a\geq0, b>0, \mu\in \mathbb{R}, \lambda>0, \alpha_1>\alpha_2>-2$, $4+2\alpha_i=2^*(\alpha_i)-2\ (i=1,2)$ with $2^*(\alpha_i)=\frac{2(N+\alpha_i)}{N-2} $ $(N=3)$ being   Hardy-Sobolev ($-2<\alpha_i<0$), Sobolev ($\alpha_i=0$) or H\'{e}non-Sobolev ($\alpha_i>0$) critical exponent of the embedding $H_{0,r}^1(B)\hookrightarrow L^p(B;|x|^{\alpha_i})$. Noting that the sign of $\mu$ gives rise to a great effect on the existence of solutions. The methods rely on Nehari manifold and the mountain pass theorem. 
 \end{abstract}
\section{Introduction and preliminaries}
In this paper we consider the Kirchhoff equation
\begin{eqnarray}
\left\{
  \begin{array}{ll}
\displaystyle-\left(a+b\int_B |\nabla u|^2dx\right)\Delta u=|x|^{\alpha_1} |u|^{4+2\alpha_1}u+\mu|x|^{\alpha_2} |u|^{4+2\alpha_2}u+\lambda h(|x|) f(u) &{\rm in}\ B,\\
u\in H_{0,r}^1(B).
  \end{array}
\right.\label{1}
\end{eqnarray}
Throughout the paper, we assume that $a\geq0, b>0, \alpha_1>\alpha_2>-2, \mu\in \mathbb{R}, \lambda>0$, $h(r):=h(|x|)$ is continuous on $(0,+\infty)$ and satisfies

$(h)$\ there exists $\beta>-2$ such that
\begin{eqnarray*}
\displaystyle 0<\liminf_{|x|\rightarrow0}\frac{h(|x|)}{|x|^\beta}
\leq\limsup_{|x|\rightarrow0}\frac{h(|x|)}{|x|^\beta}<\infty.
\end{eqnarray*}
Under the assumption $(h)$, according to \cite{2007SWW-1,2007SWW-2,2020Wang-Su-AML}, we know that the embedding
\begin{eqnarray}
H_{0,r}^1(B)\hookrightarrow L^p(B;h)\label{12}
\end{eqnarray}
is continuous for $1\leq p\leq 6+2\beta$ and $6+2\beta$ is the weighted critical exponent. Furthermore, the embedding \eqref{12} is compact for $1\leq p< 6+2\beta$. The embedding constant of \eqref{12} is denoted by $S_{hp}$, where $H_{0,r}^1(B)$ is the complete space of $C_{0,r}^\infty(B)$ with respect to the norm $\|u\|:=\left(\int_B |\nabla u|^2dx\right)^{\frac{1}{2}}$ and $C_{0,r}^\infty(B)$ contains all the radial functions in $C_{0}^\infty(B)$ which denotes the collection of
smooth functions with compact support. The weighted Lebesgue space $L^p(B;h)$ is defined by
\begin{eqnarray*}
L^p(B;h):=\left\{u: B\to \mathbb{R} \bigg|  u\ {\rm is\ measurable},\ \int_Bh(|x|)|u|^pdx<\infty\right\}\ {\rm for}\ p\geq1,
\end{eqnarray*}
with the norm $\|u\|_{L^p(B;h)}:=(\int_Bh(|x|)|u|^pdx)^{\frac{1}{p}}$.
 The Kirchhoff equation stems from the stationary analogue of the equation
\begin{eqnarray*}
\rho u_{tt}-\left(\frac{P_0}{h}+\frac{E}{2L}\int_0^L\left|\frac{\partial u}{\partial x}\right|^2dx\right)\frac{\partial^2u}{\partial x^2}=0,
\end{eqnarray*}
which was raised by Kirchhoff \cite{1883Kirchhoff} in 1883 to describe the transversal oscillations of a stretched string. For the quasilinear hyperbolic equation with boundary $\Gamma$ of class $C^\infty$ on bounded domain $\Omega\subset \mathbb{R}^N$,
\begin{eqnarray}
\left\{
  \begin{array}{ll}
\displaystyle   u_{tt}+(-1)^m\cdot M\left(\int_\Omega|\nabla^m u|^2dx\right)\cdot \Delta^m u=f\
& {\rm in}\ \Omega\times[0,T],\\
\displaystyle  u=\frac{\partial u}{\partial \nu}=\cdots=\frac{\partial^{m-1} u}{\partial \nu^{m-1}}=0\ &{\rm on}\ S=\Gamma\times[0,T],\\
\displaystyle   u=\varphi, u_t=\psi\ &{\rm for}\ t=0,x\in\Omega,
  \end{array}
\right.
\label{20}
\end{eqnarray}
where the function $M(s)$ for $s>0$ belongs to the class of $C^1$ and satisfies $M(s)\geq m_0>0$. Bernstein \cite{1940Bernstein} considered the problem \eqref{20} with the case of $m=1, N=1, f=0$. S. I. Poho\v{z}aev \cite{1975Pohozaev} studied the more general model. In 1978, J. L. Lions \cite{1978Lions-J-L} gave a "abstract" situation which contains the above situations. Let $V$ be a Hilbert space with the norm $\|\cdot\|_V$ and $V'$ be the dual space of $V$. $A$ is given satisfying
\begin{eqnarray*}
\left\{
  \begin{array}{ll}
 \displaystyle A\in L(V,V'),\ A^*=A,\\
 \displaystyle \langle Av,v\rangle=k(v)\geq \alpha \|v\|_V^2,\ \ \forall v\in V.
  \end{array}
\right.
\end{eqnarray*}
He found that there exists $u\in V$ such that
\begin{eqnarray}
\left\{
  \begin{array}{ll}
u_{tt}+M(k(u(t)))Au=f\ &{\rm in}\ V',\\
  u=\varphi, u_t=\psi\ &{\rm for}\ t=0,x\in\Omega.
  \end{array}
\right.\label{21}
\end{eqnarray}
When $A=-\Delta, M(s)=a+bs$,  \eqref{21} reduces to
\begin{eqnarray*}
\left\{
  \begin{array}{ll}
 \displaystyle u_{tt}-\left(a+b\int_\Omega |\nabla u|^2dx\right)\Delta u=f(x,u)\ &{\rm in}\ \Omega,\\
 \displaystyle u=0\ &{\rm on}\ \partial\Omega,
  \end{array}
\right.
\end{eqnarray*}
where $\Omega\subset \mathbb{R}^3$ is a bounded domain. Such a problem is usually called a nonlocal problem. The  nonlocal term $\int_\Omega |\nabla u|^2dx$ yields some difficulties
 in proving  the existence of solutions since the Fr\'{e}chet derivative of the energy functional is not weak-to-weak* continuous, which leads to that the weak limit of a Palais-Smale sequence can not be directly regarded as a weak solution of the nonlocal problem.

Motivated from the ideas of Br\'ezis and Nirenberg \cite{1983BN}, the Kirchhoff equations with Hardy critical exponent, Hardy-Sobolev critical exponent or classical Sobolev critical exponent have been studied in \cite{2019TLT,2013Figueiredo,2012WTXZ,2016LL,2013XWT,2010ACF}. For the case of weighted subcritical growth, we refer to Li and Su \cite{2015LS}. H\'{e}non \cite{1973Henon} first considered the semilinear elliptic equation ($\alpha>0$)
    \begin{eqnarray}
   \left\{
  \begin{array}{ll}
    -\Delta u= |x|^\alpha u^{q-1}, \ u>0  &{\rm in} \ B,\\
  u=0 &{\rm on} \ \partial B.
  \end{array}
  \right.\label{7}
  \end{eqnarray}
For H\'{e}non type problems \eqref{7}, Smets, Su and Willem \cite{2002SSW} first found that there exist radial and non-radial ground state solutions of \eqref{7} by  variational methods. In \cite{2020Wang-Su1,2020Wang-Su2,2020Wang-Su3}, we have considered the weighted critical equation without nonlocal term that is the case  $b=0$. Within the best of our knowledge, there is no works about Kirchhoff equation involving multiple critical exponents on the unit ball $B$, especially H\'{e}non-Sobolev critical exponent. Applying the results and ideas of \cite{2007SWW-1,2007SWW-2,2020Wang-Su-AML,2020Wang-Su1,2020Wang-Su2,2020Wang-Su3, 2013XWT}, in this paper, we establish the existence of least energy solutions of the Kirchhoff equation with multiple critical exponents, including Hardy-Sobolev, Sobolev and H\'{e}non-Sobolev critical exponents. The following tabulation sums up the cases considered in this paper:
\begin{center}
{\footnotesize
\begin{tabular}{cccccc}
\toprule[1mm]
 $a$ &$\mu$ & $\alpha_1,\alpha_2$& $f$ & $\beta$ & $\lambda$ \\
\midrule[1pt] \specialrule{0em}{2pt}{2pt}
  &  & & $(f_{1-2-3})$   &$-1<\beta\leq0$  &     \\ \specialrule{0em}{2pt}{2pt}
$a=0$ & $\mu=0$  &   $\alpha_1>-1$  & $(f_{1-2-3})$  &$\beta>0, \tau>2(2+\beta)$  &$\lambda>0$ \\ \specialrule{0em}{2pt}{2pt}
 & & & $(f_{1-2-3-4})$ & $\beta>0, \tau\leq 2(2+\beta)$  &    \\ \specialrule{0em}{2pt}{2pt}
\hline \specialrule{0em}{2pt}{2pt}
$a\geq0$ & $\mu=0$ & $\alpha_1>-1$    &$(f_{1-2-3})$   &$-$  &$\exists \lambda_{a,0}^*>0$ s.t. $\lambda>\lambda_{a,0}^*$  \\ \specialrule{0em}{2pt}{2pt}
\midrule[1pt] \specialrule{0em}{2pt}{2pt}

$a\geq0$ & $\mu>0$ & $\alpha_2>-1$    &$(f_{1-2-3})$   &$-$   &$\exists \lambda_{a,\mu}^*>0$ s.t. $\lambda>\lambda_{a,\mu}^*$  \\ \specialrule{0em}{2pt}{2pt}

\midrule[1pt] \specialrule{0em}{2pt}{2pt}
  &  &  & $(f_{5-6-7})$   &$-1<\beta\leq0$  &   \\ \specialrule{0em}{2pt}{2pt}
$a=0$ & \tabincell{c}{$\exists\mu^*<0$ s.t. $\mu^*<\mu<0$}    & $\alpha_1>-1$  & $(f_{5-6-7})$   & $\beta>0, \tau>2(2+\beta)$  &$\lambda>0$  \\ \specialrule{0em}{2pt}{2pt}
  & & & $(f_{4-5-6-7})$ & $\beta>0, \tau\leq 2(2+\beta)$  &  \\ \specialrule{0em}{2pt}{2pt}

\hline \specialrule{0em}{2pt}{2pt}
$a\geq0$ & $\mu<0$ & $\alpha_1>-1$   &$(f_{5-6-7})$  &$-$  &$\exists \bar{\lambda}_{a,\mu}>0$ s.t. $\lambda>\bar{\lambda}_{a,\mu}$ \\ \specialrule{0em}{2pt}{2pt}
\bottomrule[1mm]
\end{tabular}}
\end{center}
We mention that the additional conditions $\alpha_1>\alpha_2>-2$ and $ b>0$ are assumed, which is the basic assumptions. Throughout the paper, the symbol $(f_{1-2-3})$ indicates that $(f_1),(f_2),(f_3)$. Similarly, $(f_{1-2-3-4})$ indicates that $(f_1),(f_2),(f_3),(f_4)$.

The energy functional corresponding to the equation \eqref{1} is defined as
\begin{equation} \label{2021} \begin{array}{ll}
\Phi_{a,\mu}(u) & \displaystyle = \frac{a}{2}\int_B |\nabla u|^2dx+\frac{b}{4}\left(\int_B |\nabla u|^2dx\right)^2-\frac{1}{2^*(\alpha_1)}\int_B |x|^{\alpha_1}|u|^{2^*(\alpha_1)}dx\\[3.5mm]
& \displaystyle \ \ \ \ -\frac{\mu}{2^*(\alpha_2)}\int_B|x|^{\alpha_2}|u|^{2^*(\alpha_2)}dx-\lambda\int_Bh(|x|)F(u)dx, \end{array} \end{equation}
where $2^*(\alpha)=6+2\alpha$, $F(t):=\int_0^t f(s)ds$. Under the assumptions of the following sections, the embedding \eqref{12} implies that $\Phi_{a,\mu}$ is always of class $C^1$. The Nehari manifold of $\Phi_{a,\mu}$ is defined as
\begin{eqnarray*}
\mathcal{N}_{a,\mu}:=\{u\in H_{0,r}^1(B)\setminus\{0\}|\langle \Phi_{a,\mu}'(u),u\rangle=0\}.
\end{eqnarray*}
Let
\begin{eqnarray*}
m_{a,\mu}:=\inf_{u\in \mathcal{N}_{a,\mu}}\Phi_{a,\mu}(u),
\end{eqnarray*}
\begin{eqnarray*}
\displaystyle c_{a,\mu}:=\inf_{u\in H_{0,r}^1(B)\setminus\{0\}}\max_{t\geq0}\Phi_{a,\mu}(tu),
\end{eqnarray*}
\begin{eqnarray}
\displaystyle c_{a,\mu}^*:=\inf_{\gamma\in \Gamma}\max_{t\in[0,1]}\Phi_{a,\mu}(\gamma(t)),\label{37}
\end{eqnarray}
where
\begin{eqnarray*}
\Gamma:=\{\gamma\in C([0,1], H_{0,r}^1(B))|\gamma(0)=1,\Phi_{a,\mu}(\gamma(1))<0\}.
\end{eqnarray*}
In order to obtain the desired results, we need the inequality
\begin{eqnarray}
\int_{\mathbb{R}^3}|\nabla u|^{2} dx \geqslant S_\alpha
\left(\int_{\mathbb{R}^3}|x|^{\alpha}|u|^{6+2\alpha} dx \right)^{\frac{2}{6+2\alpha}}.\label{23}
\end{eqnarray}
The inequality \eqref{23} still holds for the bounded spherically symmetric domain and shares the same best constant. Moreover, the inequality \eqref{23} is called Sobolev inequality, Hardy-Sobolev inequality, H\'{e}non-Sobolev inequality,  and $S_\alpha$ is the best Sobolev constant, Hardy-Sobolev constant, H\'{e}non-Sobolev constant as $\alpha=0,$ $-2<\alpha<0$, $\alpha>0$, respectively, see \cite{1976Talenti,2016GR1,2020Wang-Su1,2020Wang-Su2}. In fact, $S_\alpha=(3+\alpha)\left(\frac{\omega_3}{2+\alpha}
\cdot \frac{\Gamma^2\left(\frac{3+\alpha}{2+\alpha}\right)}
{\Gamma\left(\frac{2(3+\alpha)}{2+\alpha}\right)}\right)^{\frac{2+\alpha}{3+\alpha}}$ for $\alpha>-2$ , where $\omega_3$ is the surface area of the unite sphere in $\mathbb{R}^3$. It follows from \cite{1981GS,2013GGN,1983Lieb} that the best constant $S_\alpha$ can be achieved uniquely by
\begin{eqnarray}
U_{\epsilon,\alpha}(x)=\frac{(3+\alpha)^{\frac{1}{4+2\alpha}}\epsilon^{\frac{1}{2}}}{\left(\epsilon^{2+\alpha}+|x|^{2+\alpha}\right)
^{\frac{1}{2+\alpha}}},\label{24}
\end{eqnarray}
which is the unique (up to dilations) solution of the equation
\begin{eqnarray}
\left\{\begin{array}{ll}
-\Delta u=|x|^\alpha u^{5+2\alpha},\  u>0  & \text{in}\ \mathbb{R}^3, \\
u\in D_r^{1,2}(\mathbb{R}^3),\label{25}
 \end{array}
 \right.
\end{eqnarray}
where $D_r^{1,2} (\mathbb{R}^3):=\{u\in D^{1,2}(\mathbb{R}^3)|\ u\ \mbox{is\ radial}\}$ with the norm $\|u\|_{D_r^{1,2}(\mathbb{R}^3)}=\left(\int_{\mathbb{R}^3}|\nabla u|^2dx\right)^{\frac{1}{2}}$.

In this paper, both the degenerate and non-degenerate cases are considered for $\mu=0$, $\mu>0$, $\mu<0$, respectively.
The organization of this paper is as follows. In Section \ref{sec2}, we are concerned with the existence of ground state solutions for the problem \eqref{1} with single weighted critical exponent, that is the case $\mu=0$. Section \ref{sec3} is devoted to the case of $\mu>0$. Drawing support from some conclusions of single critical exponent, the case of $\mu<0$ is studied in Section \ref{sec4}. In a forthcoming paper \cite{2021DSW-1}, we have extended the results of this paper to the nonlocal problem in $\mathbb{R}^3$.

We remark that throughout the paper, some formulas hold in the sense of subsequence. We will always adopt the convention that original sequence has been replaced by subsequence and omit the subscript.

\section{The case $\mu=0$}\label{sec2}
In this section, we focus on the problem \eqref{1} with single weighted critical exponent, i.e. $\mu=0$. We consider the cases of degeneration ($a=0$) and non-degeneration ($a>0$), respectively. Let $f$ satisfy

$(f_1)$\ $\displaystyle f\in C(\mathbb{R},\mathbb{R}), \lim_{|s|\to0}\frac{f(s)}{|s|^3}=0, \lim_{|s|\to+\infty}\frac{f(s)}{|s|^{2^*(\beta)-1}}=0;$

$(f_2)$\ there exists $\tau\in(4,6+2\beta)$ such that
\begin{eqnarray*}
0<\tau F(s)\leq f(s)s,\ \ s\in \mathbb{R}\setminus\{0\};
\end{eqnarray*}

$(f_3)$\ $\frac{f(s)}{|s|^3}$ is nondecreasing on $\mathbb{R}\setminus\{0\}$;

$(f_4)$\ $\displaystyle\lim_{|s|\to+\infty}\frac{F(s)}{|s|^{2^*(\beta)-2}}=+\infty.$

Clearly, the assumption $(f_2)$ implies that $\beta>-1$. By the assumptions $(h)$ and $(f_1)$, combining with the embedding \eqref{12}, arguing as in the proof of \cite[Lemma 2.3]{2020Wang-Su3}, we can obtain the following conclusions.
\begin{lem}\label{lem2.1}
Assume that $(h),(f_1)$ hold and that $w_n\rightharpoonup w$ in $H_{0,r}^1(B)$. Then
\begin{eqnarray}
&\displaystyle \int_Bh(|x|)f(w_n)w_ndx\to \int_Bh(|x|)f(w)wdx,&\label{92}\\
&\displaystyle \int_Bh(|x|)F(w_n)dx\to \int_Bh(|x|)F(w)dx,&\label{93}\\
&\displaystyle \int_Bh(|x|)f(w_n)\varphi dx\to \int_Bh(|x|)f(w) \varphi dx,\ \forall \varphi\in H_{0,r}^1(B).&\label{94}
\end{eqnarray}
\end{lem}
\begin{lem}\label{lem3}
Assume $\alpha_1>-1$ and $(f_{1-2-3})$. Then, for any $u\in H_{0,r}^1(B)\setminus\{0\}$, there exists a unique $t_u>0$ such that $t_u u\in \mathcal{N}_{a,0}$ and $\displaystyle\Phi_{a,0}(t_u u)=\max_{t\geq0}\Phi_{a,0}(tu)$.
\end{lem}
{\bf Proof.}\ For any fixed $u\in H_{0,r}^1(B)\setminus\{0\}$, the assumption $(f_2)$ and $\alpha_1>-1$ yield that
\begin{eqnarray*}
\Phi_{a,0}(tu)&=&\frac{at^2}{2}\|u\|^2+\frac{bt^4}{4}\|u\|^4-\frac{t^{6+2\alpha_1}}{6+2\alpha_1}\int_B|x|^{\alpha_1} |u|^{6+2\alpha_1}dx-\lambda\int_B h(|x|)F(tu)dx\\
&\leq&\frac{at^2}{2}\|u\|^2+\frac{bt^4}{4}\|u\|^4-\frac{t^{6+2\alpha_1}}{6+2\alpha_1}\int_B|x|^{\alpha_1} |u|^{6+2\alpha_1}dx\\
&\to&-\infty\ \ \ \ {\rm as}\ t\to\infty.
\end{eqnarray*}
Given $\epsilon>0$, by $(f_1)$ and the embedding \eqref{12}, we have
\begin{eqnarray*}
&&\Phi_{a,0}(tu)\\ &=&\frac{at^2}{2}\|u\|^2+\frac{bt^4}{4}\|u\|^4-\frac{t^{6+2\alpha_1}}{6+2\alpha_1}\int_B|x|^{\alpha_1} |u|^{6+2\alpha_1}dx-\lambda\int_B h(|x|)F(tu)dx\\
&\geq&\frac{bt^4}{4}\|u\|^4-\frac{t^{6+2\alpha_1}}{6+2\alpha_1}\int_B|x|^{\alpha_1} |u|^{6+2\alpha_1}dx-\lambda \int_B h(|x|) \left(\frac{\epsilon}{4} |t u|^4+C(\epsilon)|tu|^{6+2\beta}\right)dx\\
&\geq&\left(\frac{b}{4}-\frac{\epsilon \lambda  }{4S_{h 4}^2}\right)t^4\|u\|^4-\frac{t^{6+2\alpha_1}S_{\alpha_1}^{-3-\alpha_1}}{6+2\alpha_1}\|u\|^{6+2\alpha_1}
-\lambda C(\epsilon)|t|^{6+2\beta}S_{h2^*(\beta)}^{-3-\beta}\|u\|^{6+2\beta},
\end{eqnarray*}
where $C(\epsilon)$ is a positive constant depending on $\epsilon$. Choosing $\epsilon<\frac{bS_{h 4}^2}{\lambda}$, the facts that $\alpha_1>-1$ and $\beta>-1$ lead to $\Phi_{a,0}(tu)>0$ for $t>0$ small.  Thus there exists $t_u>0$ such that $\displaystyle\Phi_{a,0}(t_u u)=\max_{t\geq0} \Phi_{a,0}(tu)$ and
$\frac{d}{dt}\Phi_{a,0}(t u)|_{t=t_u}=0.$ As a consequence, $t_u u\in \mathcal{N}_{a,0}$. 
That is, 
\begin{eqnarray*}
\frac{a\|u\|^2}{(t_u)^2}+b\|u\|^4-(t_u)^{2+2{\alpha_1}}\int_B|x|^{\alpha_1} |u|^{6+2{\alpha_1}}dx-\lambda\int_B h(|x|)\frac{f(t_uu)u}{(t_u)^3}dx=0.
\end{eqnarray*}
Combining with the assumption $(f_3)$, the uniqueness of $t_u$ can be verified. \hfill$\Box$

\subsection{The subcase $a=0$}
In this subsection, we deal with the degenerative $(a=0)$ single weighted critical case and prove the following theorem.
\begin{thm}\label{thm1}
Let $\mu=0, a=0, \alpha_1>-1$ and $f$ satisfy $(f_{1-2-3})$. Then,
\begin{enumerate}
  \item[(i)] for any $\lambda>0$, the equation \eqref{1} admits a ground state solution if one of the following assumptions holds:

      {\bf (1)}\ $-1<\beta\leq0$;

      {\bf (2)}\ $\tau>2(2+\beta)$ with $\beta>0$;

      {\bf (3)}\ $\tau\leq2(2+\beta)$ with $\beta>0$ and $(f_4)$;

 \item[(ii)] there exists $\lambda_{0,0}^*>0$ such that the equation \eqref{1} has a ground state solution for any $\lambda>\lambda_{0,0}^*$.
\end{enumerate}

\end{thm}

Throughout this subsection we always suppose $\mu=0, a=0, \alpha_1>-1$ and $(f_{1-2-3})$.

\begin{lem}\label{lem1}
Any a ${\rm(PS)_c}$ sequence for $\Phi_{0,0}$ with
\begin{eqnarray*}
c<\frac{1+\alpha_1}{4(3+\alpha_1)}b^{\frac{3+\alpha_1}{1+\alpha_1}}
S_{\alpha_1}^{\frac{2(3+\alpha_1)}{1+\alpha_1}}
\end{eqnarray*}
contains a convergent subsequence.
\end{lem}
{\bf Proof.}\  Let $\{u_n\}$ be a ${\rm(PS)_c}$ sequence for $\Phi_{0,0}$, i.e.,
\begin{eqnarray}
\Phi_{0,0}(u_n)\to c,\ \ \ \Phi_{0,0}'(u_n)\to 0\ \ {\rm as}\ n\to\infty.\label{2}
\end{eqnarray}
Set $\tilde{\tau}:=\min\{\tau, 6+2\alpha_1\}$. By $(f_2)$, we deduce that
\begin{eqnarray}
\left.
  \begin{array}{ll}
c+1+C\|u_n\|&\displaystyle \geq \Phi_{0,0}(u_n)-\frac{1}{\tilde{\tau}}\langle \Phi_{0,0}'(u_n),u_n\rangle\\[3mm]
&\displaystyle =\left(\frac{b}{4}-\frac{b}{\tilde{\tau}}\right)\|u_n\|^4
+\left(\frac{1}{\tilde{\tau}}-\frac{1}{6+2\alpha_1}\right)\int_B |x|^{\alpha_1} |u_n|^{6+2\alpha_1}dx\\[4mm]
&\displaystyle\ \ \ \ \ +\lambda\int_B h(|x|)\left(\frac{1}{\tilde{\tau}}f(u_n)u_n-F(u_n)\right)dx\\[4mm]
&\displaystyle  \geq\left(\frac{b}{4}-\frac{b}{\tilde{\tau}}\right)\|u_n\|^4.
  \end{array}
\right.
\label{2.5}
\end{eqnarray}
Thus $\{u_n\}$ is bounded in $H_{0,r}^1(B)$. Then there exist $u\in H_{0,r}^1(B)$ and a renamed subsequence of $\{u_n\}$  such that
\begin{eqnarray}
\left\{
  \begin{array}{ll}
u_n\rightharpoonup u\ \ &{\rm in}\ H_{0,r}^1(B);\\
u_n\rightharpoonup u\ \ &{\rm in}\ L^{6+2\alpha_1}(B;|x|^{\alpha_1});\\
u_n(x)\to u(x)\ \ &{\rm a.e.\ in}\ B.
  \end{array}
\right.\label{11}
\end{eqnarray}
Let $v_n:=u_n-u$. To end the proof, it suffices to show $\|v_n\|\to0$. Suppose by contradiction that (up to a
subsequence)
\begin{eqnarray}
\|v_n\|\to A_\infty>0\label{A}
\end{eqnarray}
 and
\begin{eqnarray}\label{B}
\int_B |x|^{\alpha_1}|v_n|^{2^*(\alpha_1)} dx \to B_\infty.
\end{eqnarray}
Using the second limit of \eqref{2}, \eqref{2021}, Br\'ezis-Lieb Lemma \cite{1983B-Lieb} and Lemma \ref{lem2.1}, we obtain
\begin{eqnarray*}
\left.
  \begin{array}{ll}
o(1)&=\langle \Phi_{0,0}'(u_n),u\rangle\\[2mm]
&\displaystyle=b\|u_n\|^2\int_B\nabla u_n\nabla udx-\int_B |x|^{\alpha_1} |u_n|^{4+2\alpha_1}u_nudx-\lambda\int_B h(|x|)f(u_n) udx\\[3.5mm]
&\displaystyle=b\left(\|v_n\|^2+\|u\|^2\right)\|u\|^2-\int_B |x|^{\alpha_1} |u|^{6+2\alpha_1}dx-\lambda\int_B h(|x|)f(u)udx+o(1).
  \end{array}
\right.
\end{eqnarray*}
Thus
\begin{eqnarray}
0=bA_\infty^2\|u\|^2+b\|u\|^4-\int_B |x|^{\alpha_1} |u|^{6+2\alpha_1}dx-\lambda\int_B h(|x|)f(u)udx.\label{4}
\end{eqnarray}
From \eqref{2021} and \eqref{4}, we infer that
\begin{eqnarray}
\left.
  \begin{array}{ll}
\Phi_{0,0}(u)&\displaystyle=\frac{b}{4}\|u\|^4-\frac{1}{6+2\alpha_1}\int_B|x|^{\alpha_1} |u|^{6+2\alpha_1}dx-\lambda\int_B h(|x|)F(u)dx\\[3.5mm]
&\displaystyle=\left(\frac{1}{4}-\frac{1}{6+2\alpha_1}\right)\int_B|x|^{\alpha_1} |u|^{6+2\alpha_1}dx+\lambda\int_B h(|x|)\left(\frac{1}{4}f(u)u-F(u)\right)dx\\[3.5mm]
&\displaystyle\ \ \ \ -\frac{b}{4}A_\infty^2\|u\|^2\\[3.5mm]
&\displaystyle\geq-\frac{b}{4}A_\infty^2\|u\|^2.
  \end{array}
\right.\label{10}
\end{eqnarray}
On the other hand, using $\Phi_{0,0}'(u_n)\to 0$, by Br\'ezis-Lieb Lemma \cite{1983B-Lieb} and Lemma \ref{lem2.1}, we have
\begin{eqnarray*}
\left.
  \begin{array}{ll}
o(1)\|u_n\|&\displaystyle=\langle \Phi_{0,0}'(u_n),u_n\rangle\\[1mm]
&\displaystyle=b\|u_n\|^4-\int_B |x|^{\alpha_1} |u_n|^{6+2\alpha}dx-\lambda\int_B h(|x|)f(u_n) u_ndx\\[3mm]
&\displaystyle=b(\|v_n\|^2+\|u\|^2)^2-\int_B |x|^{\alpha_1} |v_n|^{6+2\alpha_1}dx
-\int_B |x|^{\alpha_1} |u|^{6+2\alpha_1}dx\\[3mm]
&\displaystyle\ \ \ \ -\lambda\int_B h(|x|)f(u) udx+o(1).
  \end{array}
\right.
\end{eqnarray*}
It follows that
\begin{eqnarray}
\left.
  \begin{array}{ll}
\displaystyle bA_\infty^4+b\|u\|^4+2bA_\infty^2\|u\|^2-B_\infty
-\int_B |x|^{\alpha_1} |u|^{6+2\alpha_1}dx-\lambda\int_B h(|x|)f(u)udx=0.
  \end{array}
\right.\label{3}
\end{eqnarray}
From \eqref{4} and \eqref{3}, we deduce that
\begin{eqnarray}
bA_\infty^4+bA_\infty^2\|u\|^2-B_\infty=0.\label{5}
\end{eqnarray}
Applying the inequality \eqref{23}, we get that
\begin{eqnarray}
B_\infty\leq S_{\alpha_1}^{-3-\alpha_1}A_\infty^{6+2\alpha_1}.\label{61}
\end{eqnarray}
Combining \eqref{61} with \eqref{5},  one concludes that
\begin{eqnarray*}
bA_\infty^4+bA_\infty^2\|u\|^2\leq S_{\alpha_1}^{-3-\alpha_1}A_\infty^{6+2\alpha_1}.
\end{eqnarray*}
It is easy to see that there exists a unique $\nu_{0,0}>0$ such that
\begin{eqnarray*}
b\nu_{0,0}^4+b\nu_{0,0}^2\|u\|^2-S_{\alpha_1}^{-3-\alpha_1}\nu_{0,0}^{6+2\alpha_1}=0\ \ \ {\rm and}\ \ \ A_\infty\geq \nu_{0,0}.
\end{eqnarray*}
Moreover
\begin{eqnarray}
\nu_{0,0}\geq\left(bS_{\alpha_1}^{3+\alpha_1}\right)^{\frac{1}{2+2\alpha_1}}.\label{26}
\end{eqnarray}
Applying again Br\'ezis-Lieb Lemma \cite{1983B-Lieb} and Lemma \ref{lem2.1}, we have
\begin{eqnarray}
\Phi_{0,0}(u_n)=\Phi_{0,0}(u)+\frac{b}{4}\|v_n\|^4+\frac{b}{2}\|v_n\|^2\|u\|^2-\frac{1}{6+2\alpha_1}\int_B |x|^{\alpha_1} |v_n|^{6+2\alpha_1}dx+o(1).\label{6}
\end{eqnarray}
Using \eqref{A} and \eqref{B}, from \eqref{5} and \eqref{6}, we deduce that
\begin{eqnarray}
\Phi_{0,0}(u)=c-\frac{b}{4}A_\infty^4-\frac{b}{2}A_\infty^2\|u\|^2+\frac{b}{6+2\alpha_1}A_\infty^4
+\frac{b}{6+2\alpha_1}A_\infty^2\|u\|^2.\label{62}
\end{eqnarray}
Putting together \eqref{62} and \eqref{26}, then
\begin{eqnarray*}
\Phi_{0,0}(u)&\leq&c-\left(\frac{1}{4}-\frac{1}{6+2\alpha_1}\right)(b\nu_{0,0}^4
+b\nu_{0,0}^2\|u\|^2)-\frac{b}{4}A_\infty^2\|u\|^2\\
&\leq&c-\left(\frac{1}{4}-\frac{1}{6+2\alpha_1}\right)b\nu_{0,0}^4-\frac{b}{4}A_\infty^2\|u\|^2\\
&\leq&c-\frac{1+\alpha_1}{4(3+\alpha_1)}b^{\frac{3+\alpha_1}{1+\alpha_1}}
S_{\alpha_1}^{\frac{6+2\alpha_1}{1+\alpha_1}}-\frac{b}{4}A_\infty^2\|u\|^2\\
&<&-\frac{b}{4}A_\infty^2\|u\|^2,
\end{eqnarray*}
which contradicts to \eqref{10}. Thus $A_\infty=0$ and the proof is complete. \hfill $\Box$

In order to estimate the mountain pass level, we first give the following estimate.
\begin{lem}\label{lem2} Assume either that one of the assumptions $\mathbf{(1)}$-$\mathbf{(3)}$ of {\rm Theorem \ref{thm1}} holds, or that $\lambda>\lambda_{0,0}^*$ for some $\lambda_{0,0}^*>0$, then there holds:
\begin{eqnarray}
c_{0,0}<\frac{1+\alpha_1}{4(3+\alpha_1)}b^{\frac{3+\alpha_1}{1+\alpha_1}}
S_{\alpha_1}^{\frac{2(3+\alpha_1)}{1+\alpha_1}}.\label{8}
\end{eqnarray}
\end{lem}
{\bf Proof.}\ From the definition of $c_{0,0}$, it is enough to check that
\begin{eqnarray*}
\max_{t>0}\Phi_{0,0}(t u_{\epsilon,\alpha_1})<\frac{1+\alpha_1}{4(3+\alpha_1)}b^{\frac{3+\alpha_1}{1+\alpha_1}}
S_{\alpha_1}^{\frac{2(3+\alpha_1)}{1+\alpha_1}},
\end{eqnarray*}
where
\begin{eqnarray}
u_{\epsilon,\alpha_1}= \varphi(x)U_{\epsilon, \alpha_1}(x), \label{uepsilon}
\end{eqnarray}
and the function $\varphi \in C_{0,r}^{\infty}(B,[0,1])$ satisfies
\begin{eqnarray*}
\left\{\begin{array}{ll}
 \varphi(|x|)\equiv 1\ \ \  \mbox{for}\ |x|\leq \frac{1}{3},\\[1em]
\varphi(|x|)\equiv0\ \ \ \mbox{for}\ |x|\geq \frac{2}{3}.
 \end{array}
 \right.
\end{eqnarray*}
By careful computations, we have
\begin{eqnarray}
&&\int_B\left|\nabla u_{\epsilon,\alpha_1}\right|^2dx
=S_{\alpha_1}^{\frac{3+\alpha_1}{2+\alpha_1}}+O\left(\epsilon\right),\label{33}\\
&&\int_B|x|^{\alpha_1} |u_{\epsilon,\alpha_1}|^{6+2\alpha_1}dx
=S_{\alpha_1}^{\frac{3+\alpha_1}{2+\alpha_1}}+O\left(\epsilon^{3+\alpha_1}\right).\label{34}
\end{eqnarray}
Applying Lemma \ref{lem3}, we conclude that the maximum of $\Phi_{0,0}(t u_{\epsilon,\alpha_1})$ for $t>0$ is attained at a unique $t_{\epsilon}$, i.e. $\displaystyle\Phi_{0,0}(t_\epsilon u_{\epsilon,\alpha_1})=\max_{t\geq0}\Phi_{0,0}(t u_{\epsilon,\alpha_1})$. It is easy to verify that
\begin{eqnarray}
0<M_1\leq t_\epsilon\leq M_2<\infty\ \ \ {\rm as}\ \epsilon>0\ {\rm small\ enough}. \label{42}
\end{eqnarray}
Furthermore,
\begin{eqnarray}
\left.
  \begin{array}{ll}
&\displaystyle \max_{t\geq0}\Phi_{0,0}(t u_{\epsilon,\alpha_1})\\[3mm]
&\displaystyle  =\frac{b t_\epsilon^4}{4}\left(\int_B |\nabla u_{\epsilon,\alpha_1}|^2dx\right)^2-\frac{t_\epsilon^{6+2\alpha_1}}{6+2\alpha_1}\int_B|x|^{\alpha_1} |u_{\epsilon,\alpha_1}|^{6+2\alpha_1}dx\\[3.5mm]
&\displaystyle\ \ \ \ \ -\lambda\int_B h(|x|)F(t_\epsilon u_{\epsilon,\alpha_1})dx\\[3.5mm]
&\displaystyle \leq\max_{t\geq0}\left\{\frac{b t^4}{4}\left(\int_B |\nabla u_{\epsilon,\alpha_1}|^2dx\right)^2-\frac{t^{6+2\alpha_1}}{6+2\alpha_1}\int_B|x|^{\alpha_1} |u_{\epsilon,\alpha_1}|^{6+2\alpha_1}dx\right\}\\[5mm]
&\displaystyle \ \ \ \ \ -\lambda\int_B h(|x|)F(t_\epsilon u_{\epsilon,\alpha_1})dx\\[3.5mm]
&\displaystyle =\frac{1+\alpha_1}{4(3+\alpha_1)}b^{\frac{3+\alpha_1}{1+\alpha_1}}
S_{\alpha_1}^{\frac{2(3+\alpha_1)}{1+\alpha_1}}+O(\epsilon)-\lambda\int_B h(|x|)F(t_\epsilon u_{\epsilon,\alpha_1})dx.
  \end{array}
\right.\label{2.17}
\end{eqnarray}
From $(h)$, there exists $C>0$ such that
\begin{eqnarray}
h(|x|)\geq C|x|^\beta\ \ {\rm as}\ |x|\leq \frac{2}{3}.\label{45}
\end{eqnarray}
For $|x|\leq \epsilon<\frac{1}{3}$, using \eqref{24} and \eqref{uepsilon}, we deduce that
\begin{eqnarray}
t_\epsilon u_{\epsilon, \alpha_1}(x)=\frac{(3+\alpha_1)^{\frac{1}{4+2\alpha_1}}t_\epsilon \epsilon^{\frac{1}{2}}}{\left(\epsilon^{2+\alpha_1}+|x|^{2+\alpha_1}\right)
^{\frac{1}{2+\alpha_1}}}\geq \frac{C(\alpha_1)M_1}{\epsilon^{\frac{1}{2}}},\label{2.18}
\end{eqnarray}
where $C(\alpha_1)=2^{-\frac{1}{2+\alpha_1}}(3+\alpha_1)^{\frac{1}{4+2\alpha_1}}$. Let $\alpha_1>-1$, $-1<\beta\leq0$ or $\tau>2(2+\beta)$ with $\beta>0$, from $(f_2)$, for any fixed $\bar{\tau}\in(\max\{2(2+\beta),4\},\tau)$ and $|s|\geq 1$, there exists $C>0$ such that
\begin{eqnarray*}
F(s)\geq C |s|^{\bar{\tau}}.
\end{eqnarray*}
Thus, combining with \eqref{2.18}, using \eqref{42} and \eqref{45}, we obtain that as $\epsilon$ small enough,
\begin{eqnarray*}
\int_B h(|x|)F(t_\epsilon u_{\epsilon,\alpha_1})dx\geq C \epsilon^{3+\beta-\frac{\bar{\tau}}{2}}\int_0^1
\frac{s^{2+\beta}}{\left(1+s^{2+\alpha_1}\right)^{\frac{\bar{\tau}}{2+\alpha_1}}}ds.
\end{eqnarray*}
The fact $\bar{\tau}>\max\{2(2+\beta),4\}$ implies
\begin{eqnarray*}
0<3+\beta-\frac{\bar{\tau}}{2}<1.
\end{eqnarray*}
Thus, using \eqref{2.17}, taking $\epsilon>0$ small enough, the estimate \eqref{8} holds for any $\lambda>0$.

Assume that $(f_4)$ holds, then for any $T>0$, there exists $s_0>0$ such that $F(s)\geq T |s|^{2(\beta+2)}$ for $|s|\geq s_0$. Together with \eqref{42} and \eqref{45}, this implies
\begin{eqnarray*}
\int_Bh(|x|)F(t_\epsilon u_{\epsilon,\alpha_1})dx
\geq C T\epsilon\int_0^1\frac{s^{3+\beta-1}}
{\left(1+s^{2+\alpha_1}\right)^{\frac{2(\beta+2)}{2+\alpha_1}}}ds.
\end{eqnarray*}
According to the arbitrariness of $T$ and \eqref{2.17}, we conclude that \eqref{8} holds for any $\lambda>0$.

For $\beta>-1$, take $v_0(x)=\varphi(|x|)|x|^{-k}$ with $0<k<\frac{1}{2}$. By the assumptions $(f_{1-2-3})$, from Lemma \ref{lem3}, there exists a unique $t_\lambda>0$ such that 
\begin{eqnarray*}
\max_{t\geq0}\Phi_{0,0}(t v_0(x))
=\frac{bt_{\lambda}^4}{4}\|v_0\|^4-\frac{t_\lambda^{6+2\alpha_1}}{6+2\alpha_1}\int_B|x|^{\alpha_1}
|v_0|^{6+2\alpha_1}dx
-\lambda \int_Bh(|x|)F(t_\lambda v_0)dx,
\end{eqnarray*}
and
\begin{eqnarray*}
b\|v_0\|^4&=&(t_\lambda)^{2+2{\alpha_1}}\int_B|x|^{\alpha_1} |v_0|^{6+2{\alpha_1}}dx+\frac{\lambda}{(t_\lambda)^3}\int_B h(|x|)f(t_\lambda v_0)v_0dx\\
&\geq&(t_\lambda)^{2+2{\alpha_1}}\int_B|x|^{\alpha_1} |v_0|^{6+2{\alpha_1}}dx.
\end{eqnarray*}
This yields
\begin{eqnarray}
t_{\lambda}\leq \left(\frac{b\|v_0\|^4}
{\int_B|x|^{\alpha_1}|v_0|^{6+2\alpha_1}dx}\right)^{\frac{1}{2+2\alpha_1}}.\label{91}
\end{eqnarray}
We claim that
\begin{eqnarray}
\lim_{\lambda\rightarrow\infty} t_\lambda=0.\label{31}
\end{eqnarray}
If not, then there exists a positive sequence $\{\lambda_i\}$ satisfying $\lambda_i\rightarrow\infty$ such that
\begin{eqnarray}
t_{\lambda_i}\to t_\infty>0\ \ \ {\rm as}\ i\to \infty.\label{27}
\end{eqnarray}
Then, using \eqref{91} and $(f_1)$, by Lebesgue Theorem, we infer that
\begin{eqnarray*}
\int_Bh(|x|)f(t_{\lambda_i} v_0)v_0dx\rightarrow \int_Bh(|x|)f(t_\infty v_0)v_0dx\ \ \ {\rm as}\ {\lambda_i}\rightarrow\infty.
\end{eqnarray*}
By $(f_2)$, for $|s|\geq 1$, there exists $C>0$ such that
\begin{eqnarray*}
f(s)s\geq \tau F(s)\geq C \tau |s|^4.
\end{eqnarray*}
Thus there exists $0<\delta_\infty<1$ such that
\begin{eqnarray*}
f(t_\infty v_0)v_0\geq \tau C(\delta_{\infty})|t_\infty|^3 |v_0|^4>0\ \ \ {\rm for}\ |x|\leq \delta_\infty.
\end{eqnarray*}
Furthermore,
\begin{eqnarray}
\int_Bh(|x|)f(t_\infty v_0)v_0dx\geq \int_{B_{\delta_\infty}}h(|x|)f(t_\infty v_0)v_0dx>0,\label{28}
\end{eqnarray}
where $B_{\delta_\infty}=\{x\in\mathbb{R}^3|\ |x|\leq\delta_\infty\}$.
Since $\frac{d\Phi_{0,0}(tv_0)}{dt}\big|_{t=t_\lambda}=0$, we deduce that
\begin{eqnarray*}
\frac{t_{\lambda_i}^3b\|v_0\|^4}{\int_Bh(|x|)f(t_{\lambda_i} v_0)v_0dx}
&\geq&{\lambda_i}.
\end{eqnarray*}
Combining with (\ref{27}) and (\ref{28}), we obtain a contradiction for $i$ large enough. Then the claim \eqref{31} is true and
\begin{eqnarray*}
\lim_{\lambda\rightarrow+\infty} \max_{t\geq0} \Phi_{0,0}(tv_0)\leq0.
\end{eqnarray*}
Thus, for $\lambda$ large enough, the desired estimate \eqref{8} holds.\hfill$\Box$\\
\\
{\bf Proof of Theorem \ref{thm1}.}\ For any fixed $u\in H_{0,r}^1(B)\setminus\{0\}$, we have
\begin{eqnarray*}
\Phi_{0,0}(u)&=&\frac{b}{4}\|u\|^4-\frac{1}{6+2\alpha_1}\int_B|x|^{\alpha_1} |u|^{6+2\alpha_1}dx-\lambda\int_B h(|x|)F(u)dx\\
&\geq&\left(\frac{b}{4}-\frac{\lambda \epsilon }{4S_{h4}^2}\right)\|u\|^4-\frac{S_{\alpha_1}^{-3-\alpha_1}}{6+2\alpha_1}\|u\|^{6+2\alpha_1}-\lambda C(\epsilon)S_{h2^*(\beta)}^{-3-\beta}\|u\|^{6+2\beta}.
\end{eqnarray*}
Taking $\epsilon<\frac{bS_{h4}^2}{\lambda }$, by the facts that $\beta>-1,\alpha>-1$, then there exists $\rho_{0,0}>0$ such that
\begin{eqnarray}
\Phi_{0,0}(u)\geq \kappa_{0,0}>0\ \ \ {\rm for}\ \|u\|=\rho_{0,0}.\label{43}
\end{eqnarray}
On the other hand, by the assumption $(f_2)$ and $\alpha_1>-1$,
\begin{eqnarray*}
\Phi_{0,0}(tu)&=&\frac{bt^4}{4}\|u\|^4-\frac{t^{6+2\alpha_1}}{6+2\alpha_1}\int_B|x|^{\alpha_1} |u|^{6+2\alpha_1}dx-\lambda\int_B h(|x|)F(tu)dx\\
&\leq&\frac{bt^4}{4}\|u\|^4-\frac{t^{6+2\alpha_1}}{6+2\alpha_1}\int_B|x|^{\alpha_1} |u|^{6+2\alpha_1}dx\\
&\to&-\infty\ \ \ \ {\rm as}\ t\to\infty.
\end{eqnarray*}
Thus there exists $\omega_{0,0}\in H_{0,r}^1(B)$ satisfying $\|\omega_{0,0}\|\geq\rho_{0,0}$ such that
\begin{eqnarray}
\Phi_{0,0}(\omega_{0,0})<0.\label{44}
\end{eqnarray}
By \eqref{43} and \eqref{44}, we can define a minimax value $c_{0,0}^*\geq\kappa_{0,0}.$
According to Lemma \ref{lem2} and the proof of \cite[Chapter 4]{1996Willem}, we deduce that
\begin{eqnarray}
c_{0,0}^*=m_{0,0}=c_{0,0}<\frac{1+\alpha_1}{4(3+\alpha_1)}b^{\frac{3+\alpha_1}{1+\alpha_1}}
S_{\alpha_1}^{\frac{6+2\alpha_1}{1+\alpha_1}}.\label{54}
\end{eqnarray}
Applying the mountain pass theorem \cite{1973AR}, there exists a ${\rm(PS)_{c_{0,0}^*}}$ sequence $\{u_n\}$ such that
\begin{eqnarray*}
\Phi_{0,0}(u_n)\to c_{0,0}^*\ \ {\rm and}\ \ \Phi_{0,0}'(u_n)\to 0\ {\rm in}\ (H_{0,r}^1(B))^{-1}\ \ {\rm as}\ n\to+\infty.
\end{eqnarray*}
Then, using \eqref{54} and Lemma \ref{lem1}, there exists $u_0\in H_{0,r}^1(B)$ such that
\begin{eqnarray*}
\Phi_{0,0}(u_n)\to \Phi_{0,0}(u_0)=m_{0,0},\ \ \Phi_{0,0}'(u_n)\to \Phi_{0,0}'(u_0)=0\ \ {\rm in}\ (H_{0,r}^1(B))^{-1}.
\end{eqnarray*}
Thus $u_0$ is a ground state solution for the equation \eqref{1}.\hfill $\Box$

\subsection{The subcase $a>0$}
In this subsection, we deal with the non-degenerative $(a>0)$ single weighted critical case and prove the following theorem.

\begin{thm}\label{thm2}
 Let $\mu=0, a>0, \alpha_1>-1$ and $f$ satisfy $(f_{1-2-3})$. There exists $\lambda_{a,0}^*>0$ such that the equation \eqref{1} has a ground state solution for any $\lambda>\lambda_{a,0}^*$.
\end{thm}

We always assume $\mu=0, a>0, \alpha_1>-1, (f_{1-2-3})$ in this subsection.

\begin{lem}\label{lem6}
Any a ${\rm (PS)_c}$ sequence of $\Phi_{a,0}$ with
\begin{eqnarray*}
c\leq\frac{2+\alpha_1}{2(3+\alpha_1)}a^{\frac{3+\alpha_1}{2+\alpha_1}}S_{\alpha_1}^{\frac{3+\alpha_1}
{2+\alpha_1}}+\frac{1+\alpha_1}{4(3+\alpha_1)}b^{\frac{3+\alpha_1}{1+\alpha_1}}
S_{\alpha_1}^{\frac{6+2\alpha_1}{1+\alpha_1}}
\end{eqnarray*}
contains a convergent subsequence.
\end{lem}
{\bf Proof.}\  Let $\{u_n\}\subset H_{0,r}^1(B)$ be a sequence such that
\begin{eqnarray}
\Phi_{a,0}(u_n)\to c,\ \ \ \Phi_{a,0}'(u_n)\to 0\ \ \ {\rm as}\ n\to\infty.\label{39}
\end{eqnarray}
Set $\tilde{\tau}:=\min\{\tau, 6+2\alpha_1\}$, then $\tilde{\tau}>4$ and
\begin{eqnarray*}
c+1+C\|u_n\|&\geq& \Phi_{a,0}(u_n)-\frac{1}{\tilde{\tau}}\langle \Phi_{a,0}'(u_n),u_n\rangle\\
&\geq&\left(\frac{a}{2}-\frac{a}{\tilde{\tau}}\right)\|u_n\|^2+
\left(\frac{b}{4}-\frac{b}{\tilde{\tau}}\right)\|u_n\|^4.
\end{eqnarray*}
Thus $\{u_n\}$ is bounded in $H_{0,r}^1(B)$. Then there exist $u\in H_{0,r}^1(B)$ and a renamed subsequence of $\{u_n\}$  such that \eqref{11} holds. Let $v_n:=u_n-u$. To end the proof, it is enough to prove $\|v_n\|\to0$. Assume for contradiction that, up to a subsequence, \eqref{A} and \eqref{B} hold. By \eqref{11}, \eqref{92}, \eqref{39} and Br\'ezis-Lieb Lemma, we have
\begin{eqnarray}
a\|u\|^2+bA_\infty^2\|u\|^2+b\|u\|^4-\int_B |x|^{\alpha_1} |u|^{6+2{\alpha_1}}dx-\lambda\int_B h(|x|)f(u)udx=0.\label{41}
\end{eqnarray}
Then we infer that
\begin{eqnarray}
\left.
  \begin{array}{ll}
\Phi_{a,0}(u) \geq-\frac{b}{4}A_\infty^2\|u\|^2.
  \end{array}
\right.\label{17}
\end{eqnarray}
On the other hand, by Br\'ezis-Lieb Lemma and $\Phi_{a,0}'(u_n)\to 0$, we have
\begin{eqnarray}
\left.
  \begin{array}{ll}
& aA_\infty^2+a\|u\|^2+bA_\infty^4+b\|u\|^4+2bA_\infty^2\|u\|^2-B_\infty\\[3mm]
&\ \ \ \displaystyle=\int_B |x|^{\alpha_1} |u|^{6+2{\alpha_1}}dx+\lambda\int_B h(|x|)f(u)udx.
  \end{array}
\right.\label{40}
\end{eqnarray}
From \eqref{41} and \eqref{40}, we deduce that
\begin{eqnarray}
aA_\infty^2+bA_\infty^4+bA_\infty^2\|u\|^2-B_\infty=0.\label{16}
\end{eqnarray}
The inequality \eqref{23} and \eqref{16} imply that
\begin{eqnarray*}
aA_\infty^2+bA_\infty^4+bA_\infty^2\|u\|^2\leq S_{\alpha_1}^{-3-\alpha_1}A_\infty^{6+2\alpha_1}.
\end{eqnarray*}
It is obvious that there exists a unique $\nu_{a,0}>0$ such that
\begin{eqnarray*}
A_\infty\geq \nu_{a,0}\ \ \ {\rm and}\ \ \ a\nu_{a,0}^2+b\nu_{a,0}^4+b\nu_{a,0}^2\|u\|^2-S_{\alpha_1}^{-3-\alpha_1}\nu_{a,0}^{6+2\alpha_1}=0.
\end{eqnarray*}
It follows that
\begin{eqnarray*}
\nu_{a,0}>(aS_{\alpha_1}^{3+\alpha_1})^{\frac{1}{4+2\alpha_1}}\ \ \ {\rm and}\ \ \
\nu_{a,0}>\left(bS_{\alpha_1}^{3+\alpha_1}\right)^{\frac{1}{2+2\alpha_1}}.
\end{eqnarray*}
Using Br\'ezis-Lieb Lemma and Lemma \ref{lem2.1}, we get that
\begin{eqnarray}
\left.
  \begin{array}{ll}
\Phi_{a,0}(u_n)&=\Phi_{a,0}(u)+\frac{a}{2}\|v_n\|^2+\frac{b}{4}\|v_n\|^4+\frac{b}{2}\|v_n\|^2\|u\|^2\\[2mm]
&\displaystyle\ \ \ \ \  -\frac{1}{6+2\alpha_1}\int_B |x|^{\alpha_1} |v_n|^{6+2\alpha_1}dx+o(1).
  \end{array}
\right.\label{15}
\end{eqnarray}
By \eqref{15} and \eqref{16}, we have
\begin{eqnarray*}
\Phi_{a,0}(u)&=&c-\frac{a}{2}A_\infty^2-\frac{b}{4}A_\infty^4-\frac{b}{2}A_\infty^2\|u\|^2
+\frac{a}{6+2\alpha_1}A_\infty^2\\
&&+\frac{b}{6+2\alpha_1}A_\infty^4+\frac{b}{6+2\alpha_1}A_\infty^2\|u\|^2.
\end{eqnarray*}
It follows that
\begin{eqnarray*}
\Phi_{a,0}(u)
&\leq&c-\left(\frac{1}{2}-\frac{1}{6+2\alpha_1}\right)a\nu_{a,0}^2-\left(\frac{1}{4}
-\frac{1}{6+2\alpha_1}\right)b\nu_{a,0}^4-\frac{b}{4}A_\infty^2\|u\|^2\\
&<&c-\frac{2+\alpha_1}{2(3+\alpha_1)}a^{\frac{3+\alpha_1}{2+\alpha_1}}
S_{\alpha_1}^{\frac{3+\alpha_1}{2+\alpha_1}}
-\frac{1+\alpha_1}{4(3+\alpha_1)}b^{\frac{3+\alpha_1}{1+\alpha_1}}
S_{\alpha_1}^{\frac{6+2\alpha_1}{1+\alpha_1}}
-\frac{b}{4}A_\infty^2\|u\|^2\\
&\leq&-\frac{b}{4}A_\infty^2\|u\|^2.
\end{eqnarray*}
This contradicts \eqref{17}. Thus $A_\infty=0$ and there exists a convergent subsequence. \hfill $\Box$

\begin{lem}\label{lem5}
There exists $\lambda_{a,0}^*>0$ such that for every $\lambda>\lambda_{a,0}^*$,
\begin{eqnarray*}
c_{a,0}^*\leq\frac{2+\alpha_1}{2(3+\alpha_1)}a^{\frac{3+\alpha_1}{2+\alpha_1}}
S_{\alpha_1}^{\frac{3+\alpha_1}{2+\alpha_1}}
+\frac{1+\alpha_1}{4(3+\alpha_1)}b^{\frac{3+\alpha_1}{1+\alpha_1}}
S_{\alpha_1}^{\frac{6+2\alpha_1}{1+\alpha_1}}.
\end{eqnarray*}
\end{lem}
{\bf Proof.}\ Arguing as what we have done in the proof of Lemma \ref{lem2}, the desired estimate on $c_{a,0}^*$ can be verified and we omit the details.\hfill $\Box$\\
\\
{\bf Proof of Theorem \ref{thm2}.}\ The proof proceeds exactly as in Theorem \ref{thm1}, the ground state solution of the equation \eqref{1} will be obtained.\hfill $\Box$

\section{The case $\mu>0$}\label{sec3}
In this section, we consider the double weighted critical Kirchhoff equation \eqref{1} with $\mu>0$. The main result of this section, which includes the cases of degeneration ($a=0$) and non-degeneration ($a>0$), is as follow.

\begin{thm}\label{thm4} Let $\mu>0, a\geq0,\alpha_2>-1$ and $f$ satisfy $(f_{1-2-3})$.
There exists $\lambda_{a,\mu}^*>0$ such that the equation \eqref{1} has a ground state solution for any $\lambda>\lambda_{a,\mu}^*$.
\end{thm}

 For the case $\mu>0$, we merely consider the case $\mu=1$. Indeed, the results still hold for the case $\mu>0$, by using the proofs as in the case $\mu=1$. Throughout this section, we always assume $\mu=1, a\geq0, \alpha_2>-1$ and $(f_{1-2-3})$. In order to establish the existence of ground state solutions, we first show the following inequalities.
\begin{lem}[{\cite[Lemmas 7.3 and 7.4]{2020Wang-Su1}}]\label{lem11}
Assume $\gamma>\xi>-2$, then for any $u\in H^1_{0,r}(B)$, we have
\begin{eqnarray*}
\|u\|_{L^{2^*(\gamma)}(B;|x|^\gamma)}\leq \tilde{C}\|u\|^{1-\theta} \|u\|_{L^{2^*(\xi)}(B;|x|^\xi)}^\theta,
\end{eqnarray*}
\begin{eqnarray*}
\|u\|_{L^{2^*(\xi)}(B;|x|^\xi)}
\leq S_\delta^{\frac{\varsigma-1}{2}}\|u\|^{1-\varsigma} \|u\|_{L^{2^*(\gamma)}(B;|x|^\gamma)}^\varsigma,
\end{eqnarray*}
where
\begin{eqnarray*}
\tilde{C}=\omega_3^{\frac{\theta-1}{2}},
\theta=\frac{2^*(\xi)}{2^*(\gamma)},
\delta=\frac{\xi 2^*(\gamma)-\gamma m}{2^*(\gamma)-m}, \varsigma=\frac{m}{2^*(\xi)}, m\in\left(0,\frac{2^*(\gamma)(2+\xi)}{2+\gamma}\right].
\end{eqnarray*}
\end{lem}

\begin{lem}\label{lem14}
Any a ${\rm(PS)_c}$ sequence of $\Phi_{a,1}$ with
\begin{eqnarray*}
c\leq \frac{2+\alpha_2}{2(3+\alpha_2)}a\bar{\nu}_{a,1}^2
+\frac{1+\alpha_2}{4(3+\alpha_2)}b\tilde{\nu}_{a,1}^4
\end{eqnarray*}
contains a convergent subsequence, where $\bar{\nu}_{a,1}$ and $\tilde{\nu}_{a,1}$ satisfy
\begin{eqnarray*}
\tilde{C}^{6+2\alpha_1}\bar{\nu}_{a,1}^{4+2\alpha_1}+\bar{\nu}_{a,1}^{4+2\alpha_2}-a
S_{\alpha_2}^{3+\alpha_2}=0,\ \ \tilde{C}^{6+2\alpha_1}\tilde{\nu}_{a,1}^{2+2\alpha_1}+\tilde{\nu}_{a,1}^{2+2\alpha_2}-b S_{\alpha_2}^{3+\alpha_2}=0.
\end{eqnarray*}
\end{lem}
{\bf Proof.}\ Let $\{u_n\}$ be a ${\rm(PS)_c}$ sequence for $\Phi_{a,1}$, namely
\begin{eqnarray}
\Phi_{a,1}(u_n)\to c,\ \ \ \Phi_{a,1}'(u_n)\to 0\ \ \ {\rm as}\ n\to\infty.\label{57}
\end{eqnarray}
As in the proof of \eqref{2.5}, the boundedness of $\{u_n\}$ in $H_{0,r}^1(B)$ follows. Thus
there exist $u\in H_{0,r}^1(B)$ and a renamed subsequence of $\{u_n\}$  such that
\begin{eqnarray}
\left\{
  \begin{array}{ll}
u_n\rightharpoonup u\ \ &{\rm in}\ H_{0,r}^1(B);\\
u_n\rightharpoonup u\ \ &{\rm in}\ L^{6+2\alpha_i}(B;|x|^{\alpha_i}),\ i=1,2;\\
u_n(x)\to u(x)\ \ &{\rm a.e.\ in}\ B.
  \end{array}
\right.\label{13}
\end{eqnarray}
Set $v_n:=u_n-u$. We claim that $\|v_n\|\to0$. Indeed, if this is false, then up to a subsequence, we have that \eqref{A}, \eqref{B} and
\begin{eqnarray}
\lim_{n\to\infty}\int_B |x|^{\alpha_2} |v_n|^{6+2\alpha_2}dx=C_\infty\label{C}
\end{eqnarray}
hold. Applying \eqref{13}, \eqref{92} and Br\'ezis-Lieb Lemma, we deduce from \eqref{57} that
\begin{eqnarray}
\displaystyle a\|u\|^2+bA_\infty^2\|u\|^2+b\|u\|^4-\sum_{i=1}^2\int_B |x|^{\alpha_i} |u|^{6+2\alpha_i}dx-\lambda\int_B h(|x|)f(u)udx=o(1).\label{55}
\end{eqnarray}
By \eqref{55}, we infer that
\begin{eqnarray}
\Phi_{a,1}(u)\geq-\frac{b}{4}A_\infty^2\|u\|^2.\label{9}
\end{eqnarray}
Combining Br\'ezis-Lieb Lemma and the fact $\Phi_{a,1}'(u_n)\to 0$ in $(H_{0,r}^1(B))^{-1}$, we have
\begin{eqnarray}
\left.
  \begin{array}{ll}
0&=\displaystyle aA_\infty^2+a\|u\|^2+bA_\infty^4+b\|u\|^4+2bA_\infty^2\|u\|^2-B_\infty-C_\infty\\
&\displaystyle\ \ \ \ \ -\sum_{i=1}^2\int_B |x|^{\alpha_i} |u|^{6+2\alpha_i}dx-\lambda\int_B h(|x|)f(u)udx.
  \end{array}
\right.\label{56}
\end{eqnarray}
By \eqref{55} and \eqref{56}, we obtain
\begin{eqnarray}
aA_\infty^2+bA_\infty^4+bA_\infty^2\|u\|^2-B_\infty-C_\infty=0.\label{58}
\end{eqnarray}
Using Lemma \ref{lem11} and \eqref{58}, we have
\begin{eqnarray*}
B_\infty>0\ \ \ {\rm and}\ \ \ C_\infty>0.
\end{eqnarray*}
Applying again Lemma \ref{lem11} and \eqref{58}, we get that
\begin{eqnarray}
\left.
  \begin{array}{ll}
bS_{\alpha_2}^2C_\infty^{\frac{2}{3+\alpha_2}}&\leq bA_\infty^4\\[1mm]
&<aA_\infty^2+bA_\infty^4+bA_\infty^2\|u\|^2\\[1mm]
&=B_\infty+C_\infty\\[1mm]
&\leq \tilde{C}^{6+2\alpha_1}A_\infty^{2(\alpha_1-\alpha_2)}
C_\infty+C_\infty
 \end{array}
\right.\label{59}
\end{eqnarray}
and
\begin{eqnarray}
\left.
  \begin{array}{ll}
aS_{\alpha_2}C_\infty^{\frac{1}{3+\alpha_2}}&\leq aA_\infty^2\\[1mm]
&<aA_\infty^2+bA_\infty^4+bA_\infty^2\|u\|^2\\[1mm]
&=B_\infty+C_\infty\\[1mm]
&\leq \tilde{C}^{6+2\alpha_1}A_\infty^{2(\alpha_1-\alpha_2)}
C_\infty+C_\infty.
  \end{array}
\right.\label{60}
\end{eqnarray}
From \eqref{59} and \eqref{60}, we have
\begin{eqnarray*}
\tilde{C}^{6+2\alpha_1}A_\infty^{2+2\alpha_1}+A_\infty^{2+2\alpha_2}-b S_{\alpha_2}^{3+\alpha_2}>0,\\ \tilde{C}^{6+2\alpha_1}A_\infty^{4+2\alpha_1}+A_\infty^{4+2\alpha_2}-a
S_{\alpha_2}^{3+\alpha_2}>0.
\end{eqnarray*}
Thus there exist $\tilde{\nu}_{a,1}:=\nu_{a,1}(\tilde{C},\alpha_1,\alpha_2,b,S_{\alpha_2})$ and $\bar{\nu}_{a,1}:=\nu_{a,1}(\tilde{C},\alpha_1,\alpha_2,a, S_{\alpha_2})$ such that
\begin{eqnarray*}
A_\infty>\tilde{\nu}_{a,1},\ \ \  \tilde{C}^{6+2\alpha_1}\tilde{\nu}_{a,1}^{2+2\alpha_1}+\tilde{\nu}_{a,1}^{2+2\alpha_2}-b S_{\alpha_2}^{3+\alpha_2}=0,
\end{eqnarray*}
and
\begin{eqnarray*}
A_\infty>\bar{\nu}_{a,1},\ \ \   \tilde{C}^{6+2\alpha_1}\bar{\nu}_{a,1}^{4+2\alpha_1}+\bar{\nu}_{a,1}^{4+2\alpha_2}-a
S_{\alpha_2}^{3+\alpha_2}=0.
\end{eqnarray*}
Using Br\'ezis-Lieb Lemma, \eqref{93} and \eqref{58}, we get
\begin{eqnarray*}
\Phi_{a,1}(u)
&=&c-\left(\frac{a}{2}-\frac{a}{6+2\alpha_2}\right)A_\infty^2
-\left(\frac{b}{4}-\frac{b}{6+2\alpha_2}\right)A_\infty^4
-\left(\frac{b}{4}-\frac{b}{6+2\alpha_2}\right)A_\infty^2\|u\|^2\\
&&-\left(\frac{1}{6+2\alpha_2}-\frac{1}{6+2\alpha_1}\right)B_\infty-\frac{b}{4}A_\infty^2\|u\|^2.
\end{eqnarray*}
It follows that
\begin{eqnarray*}
\Phi_{a,1}(u)
<c-\frac{2+\alpha_2}{2(3+\alpha_2)}a\bar{\nu}_{a,1}^2
-\frac{1+\alpha_2}{4(3+\alpha_2)}b\tilde{\nu}_{a,1}^4
-\frac{b}{4}A_\infty^2\|u\|^2
\leq-\frac{b}{4}A_\infty^2\|u\|^2.
\end{eqnarray*}
This contradicts \eqref{9}.
Thus $A_\infty=0$ and the proof is complete. \hfill $\Box$

The following lemmas are important for proving Theorem \ref{thm4}, which can be obtained by working as in Lemma \ref{lem3} and Lemma \ref{lem2}, respectively.
\begin{lem}\label{lem15}
For any $u\in H_{0,r}^1(B)\setminus\{0\}$, there exists a unique $t_u>0$ such that $t_u u\in \mathcal{N}_{a,1}$ and $\displaystyle\Phi_{a,1}(t_u u)=\max_{t\geq0} \Phi_{a,1}(tu)$.
\end{lem}

\begin{lem}\label{lem16}
There exist $\lambda_{a,1}^*>0$  and nonzero $u\in H_{0,r}^1(B)$ such that for every $\lambda>\lambda_{a,1}^*$,
\begin{eqnarray}
\max_{t\geq0}\Phi_{a,1}(tu)\leq  \frac{2+\alpha_2}{2(3+\alpha_2)}a\bar{\nu}_{a,1}^2
+\frac{1+\alpha_2}{4(3+\alpha_2)}b\tilde{\nu}_{a,1}^4,\label{69}
\end{eqnarray}
where $\bar{\nu}_{a,1}$ and $\tilde{\nu}_{a,1}$ are given in {\rm Lemma \ref{lem14}}.
\end{lem}
{\bf Proof of Theorem \ref{thm4}.}\ Using the arguments as in Theorem \ref{thm1} one can obtain the result.\hfill $\Box$

\section{The case $\mu<0$}\label{sec4}

In this section, we study the double weighted critical Kirchhoff equation \eqref{1} with $\mu<0$. The cases of degeneration ($a=0$) and non-degeneration ($a>0$) are considered, respectively. Throughout this section, we always assume $\mu<0, a\geq0, \alpha_1>-1$, and $f$ satisfies

$(f_5)$\ $f\in C(\mathbb{R},\mathbb{R}), \displaystyle\lim_{|s|\to0}\frac{f(s)}{|s|^{\zeta-1}}=0, \displaystyle\lim_{|s|\to+\infty}\frac{f(s)}{|s|^{2^*(\beta)-1}}=0$, where $\zeta:=\max\{4,6+2\alpha_2\}.$

$(f_6)$\ there exists $\tau\in(\zeta,6+2\beta)$ such that
\begin{eqnarray*}
0<\tau F(s)\leq f(s)s,\ \ s\in \mathbb{R}\setminus\{0\}.
\end{eqnarray*}

$(f_7)$\ $\frac{f(s)}{|s|^{\zeta-1}}$ is nondecreasing on $(-\infty,0)$ and $(0,+\infty)$.

 We mention that the definition of $\zeta$ in $(f_5)$ and the numerical range of $\tau$ in $(f_6)$ imply that $\beta>\max\{\alpha_2,-1\}$.
\begin{rem}\label{rem4.1}
The results of {\rm Lemma \ref{lem2.1}} still hold if the assumption $(f_1)$ is replaced by $(f_5)$.
\end{rem}

\begin{lem}\label{lem20}
Any a ${\rm(PS)_c}$ sequence of $\Phi_{a,\mu}$ with
\begin{eqnarray*}
c\leq \frac{2+\alpha_1}{2(3+\alpha_1)}a^{\frac{3+\alpha_1}{2+\alpha_1}}S_{\alpha_1}^{\frac{3+\alpha_1}
{2+\alpha_1}}+\frac{1+\alpha_1}{4(3+\alpha_1)}b^{\frac{3+\alpha_1}{1+\alpha_1}}
S_{\alpha_1}^{\frac{6+2\alpha_1}{1+\alpha_1}}
\end{eqnarray*}
contains a convergent subsequence.
\end{lem}
{\bf Proof.}\  Let $\{u_n\}\subset H_{0,r}^1(B)$ be a ${\rm (PS)_c}$ sequence for $\Phi_{a,\mu}$, that is
\begin{eqnarray}
\Phi_{a,\mu}(u_n)\to c,\ \ \ \Phi_{a,\mu}'(u_n)\to 0\ \ {\rm as}\ n\to\infty.\label{85}
\end{eqnarray}
Set $\bar{\rho}:=\min\{\rho, 6+2\alpha_1\}$, then $\bar{\rho}>4$ and
\begin{eqnarray*}
c+1+C\|u_n\|\geq \Phi_{a,\mu}(u_n)-\frac{1}{\bar{\rho}}\langle \Phi_{a,\mu}'(u_n),u_n\rangle\geq\left(\frac{b}{4}-\frac{b}{\bar{\rho}}\right)\|u_n\|^4.
\end{eqnarray*}
Thus $\{u_n\}$ is bounded in $H_{0,r}^1(B)$ and \eqref{13} holds. Let $v_n:=u_n-u$. We claim that $\|v_n\|\to0$. If the claim is not true, then we can get a subsequence such that \eqref{A}, \eqref{B} and \eqref{C} hold. By \eqref{13}, \eqref{92} and Br\'ezis-Lieb Lemma, it follows from \eqref{85} that
\begin{eqnarray}
\left.
  \begin{array}{ll}
0=&\displaystyle a\|u\|^2+bA_\infty^2\|u\|^2+b\|u\|^4-\int_B |x|^{\alpha_1} |u|^{6+2\alpha_1}dx\\[3mm]
&\displaystyle-\mu\int_B |x|^{\alpha_2} |u|^{6+2\alpha_2}dx-\lambda\int_B h(|x|)f(u)udx.
  \end{array}
\right.\label{87}
\end{eqnarray}
We deduce from \eqref{87} that
\begin{eqnarray}
\Phi_{a,\mu}(u)\geq-\frac{b}{4}A_\infty^2\|u\|^2.\label{14}
\end{eqnarray}
On the other hand, combining Br\'ezis-Lieb Lemma and $\Phi_{a,\mu}'(u_n)\to 0$ in $(H_{0,r}^1(B))^{-1}$, we have
\begin{eqnarray}
\left.
  \begin{array}{ll}
0=&\displaystyle aA_\infty^2+a\|u\|^2+b(A_\infty^2+\|u\|^2)^2-B_\infty-\int_B |x|^{\alpha_1} |u|^{6+2\alpha_1}dx\\[3mm]
&\displaystyle-\mu C_\infty-\mu\int_B |x|^{\alpha_2} |u|^{6+2\alpha_2}dx-\lambda\int_B h(|x|)f(u)udx.
  \end{array}
\right.\label{86}
\end{eqnarray}
From \eqref{87} and \eqref{86}, we obtain
\begin{eqnarray}
aA_\infty^2+bA_\infty^4+bA_\infty^2\|u\|^2-B_\infty
-\mu C_\infty=0.\label{88}
\end{eqnarray}
It follows from Lemma \ref{lem11} that
\begin{eqnarray*}
B_\infty>0\ \ \  {\rm and}\ \ \  C_\infty>0.
\end{eqnarray*}
Using \eqref{88}, we infer that
\begin{eqnarray}
bS_{\alpha_1}^2B_\infty^{\frac{2}{3+\alpha_1}}\leq bA_\infty^4<aA_\infty^2+bA_\infty^4+bA_\infty^2\|u\|^2
=B_\infty+\mu C_\infty<B_\infty,\label{89}
\end{eqnarray}
and
\begin{eqnarray}
aS_{\alpha_1}B_\infty^{\frac{1}{3+\alpha_1}}\leq aA_\infty^2<aA_\infty^2+bA_\infty^4+bA_\infty^2\|u\|^2
=B_\infty+\mu C_\infty<B_\infty,\label{90}
\end{eqnarray}
from which we deduce that
\begin{eqnarray*}
B_\infty^{\frac{1+\alpha_1}{3+\alpha_1}}>bS_{\alpha_1}^2,\ \ \
B_\infty^{\frac{2+\alpha_1}{3+\alpha_1}}>aS_{\alpha_1}.
\end{eqnarray*}
Then, by \eqref{89} and \eqref{90}, we have
\begin{eqnarray*}
&&A_\infty^4>b^{\frac{2}{1+\alpha_1}}S_{\alpha_1}^{\frac{6+2\alpha_1}{1+\alpha_1}},\ \ \  A_\infty^2>a^{\frac{1}{2+\alpha_1}}S_{\alpha_1}^{\frac{3+\alpha_1}{2+\alpha_1}}.
\end{eqnarray*}
Using again Br\'ezis-Lieb Lemma, \eqref{93} and \eqref{88}, elementary calculations show that
\begin{eqnarray*}
\Phi_{a,\mu}(u)
&=&c-\left(\frac{a}{2}-\frac{a}{6+2\alpha_1}\right)A_\infty^2
-\left(\frac{b}{4}-\frac{b}{6+2\alpha_1}\right)A_\infty^4
-\left(\frac{b}{4}-\frac{b}{6+2\alpha_1}\right)A_\infty^2\|u\|^2\\
&&-\left(\frac{\mu}{6+2\alpha_1}-\frac{\mu}{6+2\alpha_2}\right)C_\infty-\frac{b}{4}A_\infty^2\|u\|^2.
\end{eqnarray*}
Thus
\begin{eqnarray*}
\Phi_{a,\mu}(u)<c-\frac{2+\alpha_1}{2(3+\alpha_1)}a^{\frac{3+\alpha_1}{2+\alpha_1}}
S_{\alpha_1}^{\frac{3+\alpha_1}
{2+\alpha_1}}-\frac{1+\alpha_1}{4(3+\alpha_1)}b^{\frac{3+\alpha_1}{1+\alpha_1}}
S_{\alpha_1}^{\frac{6+2\alpha_1}{1+\alpha_1}}
-\frac{b}{4}A_\infty^2\|u\|^2
\leq-\frac{b}{4}A_\infty^2\|u\|^2,
\end{eqnarray*}
which contradicts \eqref{14}. The proof is complete. \hfill $\Box$

\begin{lem}\label{lem18}
For any $u\in H_{0,r}^1(B)\setminus\{0\}$, there exists a unique $t_u>0$ such that $t_u u\in \mathcal{N}_{a,\mu}$ and $\displaystyle\Phi_{a,\mu}(t_u u)=\max_{t\geq0} \Phi_{a,\mu}(tu)$.
\end{lem}
{\bf Proof.}\ The proof is similar to the proof of Lemma \ref{lem3} and we omit the details.
\hfill $\Box$

\subsection{The subcase $a=0$}
In this subsection, we deal with the degenerative case ($a=0$) and prove the following theorem.

\begin{thm}\label{thm5}Let $a=0, \alpha_1>-1$, and $f$ satisfy $(f_{5-6-7})$. Then,
\begin{enumerate}
  \item[(i)] for any $\lambda>0$, there exists $\mu^*<0$ such that for any $\mu\in(\mu^*,0)$, the equation \eqref{1} has a ground state solution under one of the following assumptions:

      {\bf (1)}\ $-1<\beta\leq0$;

      {\bf (2)}\ $\tau>2(2+\beta)$ with $\beta>0$;

      {\bf (3)}\ $\tau\leq2(2+\beta)$ with $\beta>0$ and $(f_4)$.
  \item[(ii)] for any $\mu<0$, there exists $\bar{\lambda}_{0,\mu}>0$ such that the equation \eqref{1} has a ground state solution for any $\lambda>\bar{\lambda}_{0,\mu}$.
\end{enumerate}
\end{thm}

\begin{lem}\label{lem19} Let $a=0, \alpha_1>-1$, and $f$ satisfy $(f_{5-6-7})$.
Assume either that $\mu^*<\mu<0$ for some $\mu^*<0$, one of the assumptions $(\mathbf{1})-(\mathbf{3})$ of {\rm Theorem \ref{thm5}} holds, or that $\lambda>\bar{\lambda}_{0,\mu}$ for some $\bar{\lambda}_{0,\mu}>0$, then there holds:
\begin{eqnarray}
c^*_{0,\mu}\leq\frac{1+\alpha_1}{4(3+\alpha_1)}b^{\frac{3+\alpha_1}{1+\alpha_1}}
S_{\alpha_1}^{\frac{6+2\alpha_1}{1+\alpha_1}}.\label{80}
\end{eqnarray}
\end{lem}
{\bf Proof.}\ We consider the function $u_{\epsilon,\alpha_1}$ defined in \eqref{uepsilon}, careful computations yield
\begin{eqnarray}
\int_B|x|^{\alpha_2}|u_{\epsilon,\alpha_1}|^{6+2\alpha_2}dx=\tilde{K}+O(\epsilon
^{3+\alpha_2}),\label{29}
\end{eqnarray}
where
\begin{eqnarray*}
\tilde{K}=\int_{\mathbb{R}^3}|x|^{\alpha_2}|U_{1,\alpha_1}|^{6+\alpha_2}dx.
\end{eqnarray*}
Lemma \ref{lem18} implies that
\begin{eqnarray*}
\left.
  \begin{array}{ll}
&\displaystyle
\max_{t\geq0}\Phi_{0,\mu}(tu_{\epsilon,\alpha_1})\\[3.5mm]
&\displaystyle=\frac{t_{\epsilon,\mu,\lambda}^4}{4}\|u_{\epsilon,\alpha_1}\|^4
-\frac{t_{\epsilon,\mu,\lambda}^{6+2\alpha_1}}{6+2\alpha_1}
\int_B|x|^{\alpha_1}|u_{\epsilon,\alpha_1}|^{6+2\alpha_1}dx-\frac{\mu t_{\epsilon,\mu,\lambda}^{6+2\alpha_2}}{6+2\alpha_2}
\int_B|x|^{\alpha_2}|u_{\epsilon,\alpha_1}|^{6+2\alpha_2}dx\\[3.5mm]
&\displaystyle\ \ \ \ \
-\lambda\int_Bh(|x|)F(t_{\epsilon,\mu,\lambda}u_{\epsilon,\alpha_1})dx\\[3.5mm]
&\displaystyle\leq\max_{t\geq0}\left\{\frac{t^4}{4}\|u\|^4-\frac{t^{6+2\alpha_1}}{6+2\alpha_1}
\int_B|x|^{\alpha_1}|u_{\epsilon,\alpha_1}|^{6+2\alpha_1}dx-\lambda\int_Bh(|x|)
F(tu_{\epsilon,\alpha_1})dx\right\}\\[3.5mm]
&\displaystyle \ \ \ \ \ -\frac{\mu t_{\epsilon,\mu,\lambda}^{6+2\alpha_2}}{6+2\alpha_2}
\int_B|x|^{\alpha_2}|u_{\epsilon,\alpha_1}|^{6+2\alpha_2}dx.
\end{array}
\right.
\end{eqnarray*}
It is easy to verify that $t_{\epsilon,\mu,\lambda}$ is bounded as $\epsilon\to0^+,\mu\to0^-$, and $t_{\epsilon,\mu,\lambda}$ is away from 0 as $\epsilon\to0^+$. By Lemma \ref{lem2}, we can choose $\epsilon>0$ small enough such that
\begin{eqnarray}
\left.
  \begin{array}{ll}
&\displaystyle\max_{t\geq0}\left\{\frac{t^4}{4}\|u\|^4\!-\!\frac{t^{6+2\alpha_1}}{6+2\alpha_1}
\int_B|x|^{\alpha_1}|u_{\epsilon,\alpha_1}|^{6+2\alpha_1}dx\!-\!
\lambda\int_Bh(|x|)F(tu_{\epsilon,\alpha_1})
dx\right\}\!\\[3.5mm]
&\ \ \ \displaystyle<\!\frac{1+\alpha_1}{4(3+\alpha_1)}b^{\frac{3+\alpha_1}{1+\alpha_1}}
S_{\alpha_1}^{\frac{6+2\alpha_1}{1+\alpha_1}}.
\end{array}
\right.\label{84}
\end{eqnarray}
It follows from \eqref{29} that
\begin{eqnarray*}
\frac{\mu t_{\epsilon,\mu,\lambda}^{6+2\alpha_2}}{6+2\alpha_2}
\int_B|x|^{\alpha_2}|u_{\epsilon,\alpha_1}|^{6+2\alpha_2}dx\to0\ \ \ {\rm as}\ \epsilon\to0^+,\mu\to0^-.
\end{eqnarray*}
Thus there exists $\mu^*<0$ such that \eqref{80} holds for $\mu^*<\mu<0$.

For the case of large $\lambda$, the proof is similar to that of Lemma \ref{lem2}. \hfill$\Box$\\
\\
{\bf Proof of Theorem \ref{thm5}.}\ It is easy to see that the functional $\Phi_{0,\mu}$ exhibits a mountain pass geometry. According to the mountain pass theorem \cite{1973AR}, Lemma \ref{lem19} and the proof of \cite[Chapter 4]{1996Willem}, we can define a minimax value $0<c^*_{0,\mu}\leq \frac{2+\alpha}{4(6+\alpha_1)}
b^{\frac{3+\alpha_1}{1+\alpha_1}}S_{\alpha_1}^{\frac{6+2\alpha_1}{1+\alpha_1}}$ and there exists a ${\rm(PS)_{c^*_{0,\mu}}}$ sequence $\{u_n\}$ such that
\begin{eqnarray*}
\Phi_{0,\mu}(u_n)\to c^*_{0,\mu}=c_{0,\mu}=m_{0,\mu}\ \ {\rm and}\ \ \Phi_{0,\mu}'(u_n)\to 0\ {\rm in}\ (H_{0,r}^1(B))^{-1}\ \ {\rm as}\ n\to+\infty.
\end{eqnarray*}
From Lemma \ref{lem20}, there exists a convergent subsequence such that
\begin{eqnarray*}
\Phi_{0,\mu}(u_n)\to \Phi_{0,\mu}(u_0)=m_{0,\mu},\ \ \Phi_{0,\mu}'(u_n)\to \Phi_{0,\mu}'(u_0)=0\ \ {\rm in}\ (H_{0,r}^1(B))^{-1}.
\end{eqnarray*}
Thus $u_0$ is a ground state solution for the equation \eqref{1}.\hfill $\Box$

\subsection{The subcase $a>0$}

In this subsection, we study the non-degenerative case ($a>0$) and prove the following theorem.

\begin{thm}\label{thm6} Let $a>0, \alpha_1>-1$, and $f$ satisfy $(f_{5-6-7})$.
Then there exists $\bar{\lambda}_{a,\mu}>0$ such that the equation \eqref{1} has a ground state solution for any $\lambda>\bar{\lambda}_{a,\mu}$.
\end{thm}

\begin{lem}\label{lem22} Let $a>0, \alpha_1>-1$, and $f$ satisfy $(f_{5-6-7})$.
Then there exist $\bar{\lambda}_{a,\mu}>0$ and $u\in H_{0,r}^1(B)$ such that for any $\lambda>\bar{\lambda}_{a,\mu}$,
\begin{eqnarray}
\max_{t\geq0}\Phi_{a,\mu}(tu)\leq \frac{2+\alpha_1}{2(3+\alpha_1)}a^{\frac{3+\alpha_1}{2+\alpha_1}}S_{\alpha_1}^{\frac{3+\alpha_1}
{2+\alpha_1}}+\frac{1+\alpha_1}{4(3+\alpha_1)}b^{\frac{3+\alpha_1}{1+\alpha_1}}
S_{\alpha_1}^{\frac{6+2\alpha_1}{1+\alpha_1}}.\label{96}
\end{eqnarray}
\end{lem}
{\bf Proof.}\ The proof works similarly as in Lemma \ref{lem2}, so we omit the proof. \hfill$\Box$\\

Using the same argument as in Theorem \ref{thm5}, Theorem \ref{thm6} is a consequence of Lemma \ref{lem20}, Lemma \ref{lem18} and Lemma \ref{lem22}.
\hfill$\Box$

\

\noindent{\bf Acknowledgment} \ This work is supported by NSFC12301144, NSFC12271373, NSFC12171326 and 2024NSFSC1342.

\end{document}